\documentclass{amsart}

\address{Department of Mathematics, Columbia University, New York, NY   10027, USA}
\email{sucharit.sarkar@gmail.com}

\usepackage{epsfig, color,
bm, chngcntr, graphics, psfrag}
\usepackage[all,cmtip]{xy}

\newtheorem{thm}{Theorem}[section]
\newtheorem{lem}[thm]{Lemma}

\newtheorem{defn}[thm]{Definition}

\counterwithin{figure}{section}

\newcommand{\R}{\mathbb{R}}
\newcommand{\Z}{\mathbb{Z}}
\newcommand{\F}{\mathbb{F}}
\newcommand{\C}{\mathbb{C}}

\newcommand{\wh}{\widehat}
\newcommand{\del}{\partial}
\newcommand{\wt}{\widetilde}
\newcommand{\mc}{\mathcal}
\newcommand{\mb}{\mathbb}
\newcommand{\mf}{\mathfrak}

\newcommand{\sm}{\setminus}
\newcommand{\al}{\alpha}
\newcommand{\be}{\beta}

\newcommand{\Ozsvath}{Ozsv\'{a}th}
\newcommand{\Zoltan}{Zolt\'{a}n}
\newcommand{\Szabo}{Szab\'{o}}

\newcommand{\ith}{^{\text{th}}}
\newcommand{\Hom}{\operatorname{Hom}}
\newcommand{\SpinC}{\text{Spin}^{\text{C}}}
\newcommand{\Sym}{\text{Sym}}
\newcommand{\whcfk}{\widehat{\mathit{CFK}}}
\newcommand{\whhfk}{\widehat{\mathit{HFK}}}
\newcommand{\whcfl}{\widehat{\mathit{CFL}}}
\newcommand{\whhfl}{\widehat{\mathit{HFL}}}
\newcommand{\whhfg}{\widehat{\mathit{HFG}}}

\title{A note on sign conventions in link Floer homology} 

\author{Sucharit Sarkar} 

\begin{document}

\begin{abstract}
  For knots in $S^3$, the bi-graded hat version of knot Floer homology   is defined over $\Z$; however, for an $l$-component link $L$ in   $S^3$ with $l>1$, there are $2^{l-1}$ bi-graded hat versions of link   Floer homology defined over $\Z$; the multi-graded hat version of   link Floer homology, defined from holomorphic considerations, is   only defined over $\F_2$; and there is a multi-graded version of   link Floer homology defined over $\Z$ using grid diagrams. In this   short note, we try to address this issue, by extending the   $\F_2$-valued multi-graded link Floer homology theory to $2^{l-1}$   $\Z$-valued theories. A grid diagram representing a link gives rise   to a chain complex over $\F_2$, whose homology is related to the   multi-graded hat version of link Floer homology of that link over   $\F_2$. A sign refinement of the chain complex exists, and for   knots, we establish that the sign refinement does indeed correspond   to the sign assignment for the hat version of the knot Floer   homology. For links, we create $2^{l-1}$ sign assignments on the   grid diagrams, and show that they are related to the $2^{l-1}$   multi-graded hat versions of link Floer homology over $\Z$, and one   of them corresponds to the existing sign refinement of the grid   chain complex.
\end{abstract}



\maketitle

\section{Introduction}\label{sec:introduction}

Knot Floer homology, primarily as an invariant for knots and links
inside $S^3$, was discovered by Peter \Ozsvath{} and \Zoltan{}
\Szabo{} \cite{POZSzknotinvariants}, and independently by Jacob
Rasmussen \cite{JR}. Later, a related invariant for links, called link
Floer homology, was constructed by Peter \Ozsvath{} and \Zoltan{}
\Szabo{} \cite{POZSzlinkinvariants}. However, due to certain
orientation issues, the link invariant was only constructed over
$\F_2$, instead of $\Z$. This short note is the author's effort to
understand the orientation issues that are known, and to resolve some
of the issues that are unknown.

Let us describe the algebraic structure of knot Floer homology in the
simplest case, as described in \cite{POZSzknotinvariants}. Let $K$ be
a null-homologous knot in $\#^{l-1}(S^1\times S^2)$. Then there are
$2^{l-1}$ bi-graded chain complexes over $\Z$, such that they all give
rise to the same complex when tensored with $\F_2$. The two
gradings are called the Maslov grading $M$ and the Alexander grading
$A$. The boundary maps preserve the Alexander grading, but lower the
Maslov grading by one. Therefore, the Maslov grading acts as the
homological grading while the Alexander grading acts as an extra
filtration.  The homology of the chain complexes is called the hat
version of the knot Floer homology. Therefore, we get an $\F_2$-valued
bi-graded hat version of knot Floer homology and $2^{l-1}$ $\Z$-valued
bi-graded hat versions of knot Floer homology.

The reason for working with null-homologous knots in connected sums of
$S^1\times S^2$ is very simple. We want to work with links in $S^3$.
However, a link with $l$ components in $S^3$ very naturally gives rise
to a null-homologous knot in $\#^{l-1}(S^1\times S^2)$,
\cite{POZSzknotinvariants}. Therefore, what we have is the following.
Given a link $L\subset S^3$, with $l$ components, and after making
certain auxiliary choices, we get $2^{l-1}$ bi-graded chain complexes
over $\Z$, henceforth denoted by $\whcfk(L,\Z,\mf{o})$, where
$\mf{o}$, called an orientation system, takes values in an indexing
set of $2^{l-1}$ elements, and records which of the $2^{l-1}$ chain
complexes is the one under consideration. All of the $2^{l-1}$ chain
complexes give rise the same bi-graded chain complex over $\F_2$,
$\whcfk(L,\F_2)=\whcfk(L,\Z,\mf{o})\otimes \F_2$. The reader
should be warned that these bi-graded chain complexes,
$\whcfk(L,\Z,\mf{o})$ and $\whcfk(L,\F_2)$, are not
link-invariants (they might depend on the auxiliary choices that we
did not specify, but simply alluded to), but their homologies are link
invariants. Therefore, we get one $\F_2$-valued bi-graded hat version
of knot Floer homology $\whhfk(L,\F_2)=H_*(\whcfk(L,\F_2))$, and
$2^{l-1}$ $\Z$-valued bi-graded hat versions of knot Floer homology
$\whhfk(L, \Z, \mf{o})=H_*(\whcfk(L,\Z,\mf{o}))$. We often
let $\whhfk(L,\Z)$ denote any one of the $2^{l-1}$ versions, or a canonical one, namely the one coming from the canonical choice of orientation systems in \cite{POZSzapplications}. However, to decide which of the $2^{l-1}$ groups $\whhfk(L,\Z,\mf{o})$ is the canonical one, one needs to understand some of the other versions of link Floer homology, most notably the infinity version. This seems to be a harder problem, for reasons that we will discuss shortly.

In \cite{POZSzlinkinvariants}, the story for links is treated
in a slightly different light, and a new definition of link Floer
homology is given. Given a link $L$ with $l$ components in $S^3$,
modulo certain choices, a chain complex $\whcfl(L,\F_2)$ over $\F_2$
is constructed. The chain complex carries $(l+1)$ gradings: a single
Maslov grading $M$, which is lowered by one by the boundary map, and
$l$ Alexander gradings $A_1, A_2, \ldots, A_l$, one for each link
component, each of which is preserved by the boundary map. The
homology of the chain complex $\whhfl(L,\F_2)=H_*(\whcfl(L,\F_2))$
is an $\F_2$-valued $(l+1)$-graded homology theory, called the link
Floer homology, and it is a link invariant. These two definitions,
\emph{a priori}, are different. Therefore, we have been careful
throughout; we have called the definition from
\cite{POZSzknotinvariants} the knot Floer homology (even when talking
about links), and denoted it by $\whhfk$, and we have called the
definition from \cite{POZSzlinkinvariants} the link Floer homology,
and denoted it by $\whhfl$. However, by a miraculous coincidence, it
turns out that if we condense the $l$ Alexander gradings in
$\whhfl(L,\F_2)$ into one single Alexander grading $A=\sum_i A_i$,
then the resulting $\F_2$-valued bi-graded homology group is
isomorphic to $\whhfk(L,\F_2)$.

In this note, we will complete the picture by constructing $2^{l-1}$
$\Z$-valued chain complexes, $\whcfl(L,\Z,\mf{o})$, each carrying a
Maslov grading $M$, and $l$ Alexander gradings $A_1,A_2,\ldots,A_l$,
such that the homologies
$\whhfl(L,\Z,\mf{o})=H_*(\whcfl(L,\Z,\mf{o}))$ are link
invariants, and on condensing the $l$ Alexander gradings into one
Alexander grading $A=\sum_i A_i$, we get the $2^{l-1}$ $\Z$-valued
bi-graded homology groups $\whhfk(L,\Z,\mf{o})$.

A similar story (except possibly the last bit of coincidence) holds
for the other versions of link Floer homologies, most notably the
minus, plus and infinity versions; however, the holomorphic
considerations and the orientation issues are significantly more
subtle. In particular, we will encounter boundary degenerations, and we will have to orient the relevant moduli spaces in a consistent fashion. We plan to address this problem in future work. Understanding the orientation issues for all versions of link Floer homology will help us understand which of the $2^{l-1}$ link Floer homology groups is the canonical one and whether it has some sort of a useful characterization.

For the second part of the discourse, we concentrate on the
computational aspects of the theory. Ever since knot Floer homology
saw the light of day \cite{POZSzknotinvariants}, \cite{JR},
  \cite{POZSzlinkinvariants}, and some of its immense strengths were
discovered \cite{POZSzgenusbounds}, \cite{POZSzthurstonnorm}, \cite{YN}, people were
interested in algorithms to compute it. There have been several
recent developments towards computing various versions of link Floer
homology for links in $S^3$ \cite{CMPOSS}, \cite{SSJW}, \cite{POZSzcube}, \cite{POASZSz}.
We choose to concentrate on the algorithm from \cite{CMPOSS}: the link
$L$ in $S^3$ is represented by a toroidal grid diagram $G$, such that
the $i\ith$ component is represented by $m_i$ vertical line segments
and $m_i$ horizontal line segments; an $\F_2$-valued $(l+1)$-graded
chain complex $C(G)$ is constructed such that its homology $H_*(C(G))$
is isomorphic to
$\whhfl(L,\F_2)\otimes_i(\otimes^{m_i-1}(\F_2\oplus\F_2))$, where,
in the $(\F_2\oplus\F_2)$ that is tensored with itself $(m_i-1)$
times, for one of the generators, all the $(l+1)$ gradings are zero,
and for the other generator, the Maslov grading $M=-1$, and the
Alexander grading $A_j=-\delta_{ij}$.

Very shortly thereafter, \cite{CMPOZSzDT} assigned signs of $\pm 1$ to each of the boundary maps in the chain complex $C(G)$ in a well defined way, such that it remains a chain complex and its  homology (over $\Z$) is isomorphic to $\whhfg(L,\Z)\otimes_i(\otimes^{m_i-1}(\Z\oplus\Z))$, for some $(l+1)$-graded group $\whhfg(L,\Z)$, which is a link invariant. 
A very natural question that arises is whether the new homology group $\whhfg(L,\Z)$ is isomorphic to $\whhfl(L,\Z,\mf{o})$ for some $\mf{o}$. We establish that the answer is in the affirmative, and indeed, we construct $2^{l-1}-1$ other sign assignments on the boundary maps of $C(G)$, such that the homologies of these $2^{l-1}$ sign refined grid chain complexes correspond precisely to the $2^{l-1}$ $\Z$-valued $(l+1)$-graded homology groups $\whhfl(L,\Z)$. Once again, it is an interesting question whether $\whhfg(L,\Z)$ is isomorphic to the canonical $\whhfl(L,\Z)$, and once again, we are unable to answer it with our present methods. It is also an interesting endeavor to find two $l$-component links $L_1$ and $L_2$, such that $\whcfl(L_1,\F_2)$ is isomorphic to $\whcfl(L_2,\F_2)$ as $(l+1)$-graded $\F_2$-modules, there is a bijection between the set of  $2^{l-1}$ groups $\whcfk(L_1,\Z)$ and the set of $2^{l-1}$ groups $\whcfk(L_2,\Z)$ such that the corresponding groups are isomorphic as bi-graded $\Z$-modules, $\whhfg(L_1,\Z)$ is isomorphic to $\whhfg(L_2,\Z)$ as $(l+1)$-graded $\Z$-modules, but there is no bijection between the set of $2^{l-1}$ groups $\whhfl(L_1,\Z)$ and the set of $2^{l-1}$ groups $\whhfl(L_2,\Z)$ such that the corresponding groups are isomorphic as $(l+1)$-graded $\Z$-modules.

This is a rather short paper. We expect the reader to be already
familiar with most of \cite{CMPOZSzDT}, \cite{POZSzknotinvariants},
  \cite{POZSzlinkinvariants}. Despite trying our level best to be as
self-contained as possible, we will still be rather fast in our
exposition.

\subsection*{Acknowledgment}
The work was done when the author was supported by the Clay Research Fellowship. He would like to thank Robert Lipshitz, Peter \Ozsvath{} and \Zoltan{} \Szabo{} for several helpful discussions. He would also like to thank the referee for providing useful comments and for pointing out the errors.

\section{Floer homology}\label{sec:floerhomology}

For the first part of the section, in the following few numbered paragraphs, we will briefly review the basics of Heegaard Floer homology. The interested reader is referred to \cite{POZSz}, \cite{POZSzapplications} for more details.

\subsection{} A \emph{Heegaard diagram} is an object $\mc{H}=(\Sigma_g, \alpha_1, \ldots,\alpha_{g+k-1},\beta_1,\ldots,\beta_{g+k-1},\allowbreak X_1,\ldots,X_k,O_1,\ldots,O_k)$, where: $\Sigma_g$ is a Riemann surface of genus $g$; $\alpha=(\alpha_1,\ldots,\allowbreak \alpha_{g+k-1})$ is $(g+k-1)$-tuple of disjoint simple closed curves such that $\Sigma_g\sm\alpha$ has $k$ components; $\beta=(\beta_1,\ldots, \beta_{g+k-1})$ is $(g+k-1)$-tuple of disjoint simple closed curves such that $\Sigma_g\sm\beta$ has $k$ components; the $\alpha$ circles are transverse to the $\beta$ circles; $X=(X_1,\ldots,X_k)$ is a $k$-tuple of points such that each component of $\Sigma_g\sm\alpha$ has an $X$ marking, and each component of $\Sigma_g\sm\beta$ has an $X$ marking; $O=(O_1,\ldots,O_k)$ is a $k$-tuple of points such that each component of $\Sigma_g\sm\alpha$ has an $O$ marking, and each component of $\Sigma_g\sm\beta$ has an $O$ marking; and the diagram is assumed to be \emph{admissible}, which is a technical condition that we will describe later.

\subsection{}A Heegaard diagram represents an oriented link $L$ inside a three-manifold $Y$ in the following way: the pair $(\Sigma_g,\alpha)$ represents genus $g$ handlebody $U_{\alpha}$; the pair $(\Sigma_g,\beta)$ represents genus $g$ handlebody $U_{\beta}$; the ambient three-manifold $Y$ is obtained by gluing $U_{\alpha}$ to $U_{\beta}$ along $\Sigma_g$; the $X$ markings are joined to the $O$ markings by $k$ simple oriented arcs in the complement of the $\alpha$ circles, and the interiors of the $k$ arcs are pushed slightly inside the handlebody $U_{\alpha}$; the $O$ markings are joined to the $X$ markings by $k$ simple oriented arcs in the complement of the $\beta$ circles, and the interiors of the $k$ arcs are pushed slightly inside the handlebody $U_{\beta}$; the union of these $2k$ oriented arcs is the oriented link $L$. Let the link have $l$ components, and let $2 m_i$ be the number of arcs that represent $L_i$, the $i\ith$ component of the link $L$.  Therefore, $k=\sum_i m_i \geq l$. In \cite{POZSzlinkinvariants}, the case $k=l$ is studied, and in \cite{POZSzknotinvariants}, the subcase $k=l=1$ is dealt with. We will always assume that $L_i$ is null-homologous in $Y$, for each $i$.

\subsection{}Consider $(g+k-1)$-tuples of points $x=(x_1,\ldots,x_{g+k-1})$, such that each $\alpha$ circle contains some $x_i$, and each $\beta$ circle contains some $x_j$. To each such tuple $x$, we can associate a $\SpinC$ structure $\mf{s}_x$ on the ambient three-manifold $Y$. In all the three-manifolds that we will consider, we will be interested in a canonical torsion $\SpinC$ structure.  In particular, for $Y=\#^n S^1\times S^2$, we will be interested in the unique torsion $\SpinC$ structure. A \emph{generator} is a $(g+k-1)$-tuple $x$ of the type described above, such that $\mf{s}_x$ is the canonical $\SpinC$ structure. The set of all generators in a Heegaard diagram $\mc{H}$ is denoted by $\mc{G}_{\mc{H}}$. An \emph{elementary domain} is a component of $\Sigma_g\sm(\alpha\cup\beta)$. A \emph{domain} $D$ joining a generator $x$ to a generator $y$, is a $2$-chain generated by elementary domains such that $\del(\del D|_{\alpha})=y-x$.  The set of all domains joining $x$ to $y$ is denoted by $\mc{D}(x,y)$. A \emph{periodic domain} $P$ is a $2$-chain generated by elementary domains such that $\del(\del P|_{\alpha})=0$.  The set of periodic domains is denoted by $\mc{P}_{\mc{H}}$, and there is a natural bijection between $\mc{P}_{\mc{H}}$ and $\mc{D}(x,x)$ for any generator $x$. If $D$ is a domain, and if $p$ is a point lying in an elementary domain, then $n_p(D)$ denotes the coefficient of the $2$-chain $D$ at that elementary domain. Let $n_X(D)=\sum_i n_{X_i}(D)$ and $n_O(D)=\sum n_{O_i}(D)$. Furthermore, let $n_{X,i}(D)$ denote the sum of $n_{X_j}(D)$ for all the $X_j$ markings that lie in $L_i$, and let $n_{O,i}(D)$ denote the sum of $n_{O_j}(D)$ for all the $O_j$ markings that lie in $L_i$.  A domain is said to be \emph{non-negative} if it has non-negative coefficients in every elementary domain. A domain $D$ is said to be \emph{empty} if $n_{X_i}(D)=n_{O_i}(D)=0$ for all $i$. A Heegaard diagram is called \emph{admissible} if there are no non-negative, non-trivial empty periodic domains. The set of all empty domains in $\mc{D}(x,y)$ is denoted by $\mc{D}^0(x,y)$, and the set of all empty periodic domains is denoted by $\mc{P}^0_{\mc{H}}$. The set $\mc{P}^0_{\mc{H}}$ forms a free abelian group of rank $b_1(Y)+l-1$.

\subsection{}Every domain $D$ has an integer valued \emph{Maslov index} $\mu(D)$ associated to it, which satisfies certain properties that we will mention as we need them. In all the Heegaard diagrams that we will consider, the following additional restrictions will hold: if $P\in\mc{D}(x,x)$, then $\mu(P)=2 n_O(P)$ and, since $L_i$ is null-homologous in $Y$, $n_{X,i}(P)= n_{O,i}(P)$ for all $i$. This allows us to define $(l+1)$ relative gradings. Given two generators $x,y$, choose a domain $D\in\mc{D}(x,y)$ (since $\mf{s}_x=\mf{s}_y$, the set $\mc{D}(x,y)$ is non-empty), and let the \emph{relative Maslov grading} $M(x,y)=\mu(D)-2n_O(D)$, and let the \emph{relative Alexander grading} $A_i(x,y)=n_{X,i}(D)-n_{O,i}(D)$. In certain situations, with certain additional hypotheses, these gradings can be lifted to absolute gradings. However, for convenience, we will not work with absolute gradings right away. Therefore, until Lemma \ref{lem:absolutegrading}, whenever we talk about the Maslov grading $M$, or the Alexander grading $A_i$, we mean some affine lift of the corresponding relative grading, which is only well-defined up to a translation by $\Z$. Let $Q_i=\Z\oplus\Z$ be the $(l+1)$-graded group, with the two generators lying in gradings $(0,0,\ldots,0)$ and $(-1,-\delta_{i1},\ldots,-\delta_{il})$, where $\delta$ is the Kronecker delta function.

\subsection{}For the analytical aspects of the theory, which we are about to describe now, the reader is strongly advised to read Section 3 of \cite{POZSz}.  Let $\Sym^{g+k-1}(\Sigma_g)$ be $(g+k-1)$-fold symmetric product, and let $J_s$ be a path of nearly symmetric almost complex structures on it, obtained as a small perturbation of the constant path of nearly symmetric almost complex structure $\Sym^{g+k-1}(\mf{j})$, where $\mf{j}$ is a fixed complex structure on $\Sigma_g$, such that $J_s$ achieves certain transversality that we will describe later.  The subspaces $\mb{T}_{\alpha}=\alpha_1\times\cdots\times\alpha_{g+k-1}$ and $\mb{T}_{\beta}=\beta_1\times\cdots\times\beta_{g+k-1}$ are two totally real tori. Notice that $\mc{G}_{\mc{H}}$ is in a natural bijection with a subset of $\mb{T}_{\alpha}\cap\mb{T}_{\beta}$. Fix $\mf{p}>2$. Given a domain $D\in\mc{D}(x,y)$, let $\mc{B}(D)$ be the space of all $L^{\mf{p}}_1$ maps $u$ from $[0,1]\times\R\subset\C$ to $\Sym^{g+k-1}(\Sigma_g)$, such that: $u$ maps $\{0\}\times\R$ to $\mb{T}_{\alpha}$; $u$ maps $\{1\}\times\R$ to $\mb{T}_{\beta}$; $\lim_{t\rightarrow\infty}u(s+it)=x$ with a certain pre-determined asymptotic behavior; $\lim_{t\rightarrow -\infty}u(s+it)=y$ with a certain pre-determined asymptotic behavior; for any point $p$ in any elementary domain, the algebraic intersection number between $u$ and $\{p\}\times \Sym^{g+k-2}(\Sigma_g)$ is $n_p(D)$, or, as it is colloquially stated, the domain $D$ is the \emph{shadow} of $u$. \Ozsvath{} and \Szabo{} define a vector bundle $\mc{L}$ over $\mc{B}(D)$, and a section $\xi$ of that bundle depending on $J_s$, such that the linearization of the section $D_u\xi$ is a Fredholm operator for every $u\in\mc{B}(D)$. The transversality of the path $J_s$ that we mentioned earlier, simply means that the Fredholm section $\xi$ is transverse to the $0$-section of $\mc{L}$.  The intersection of $\xi$ and the $0$-section is denoted by $\mc{M}_{J_s}(D)$, and it consists precisely of the $J_s$-holomorphic maps. There is an $\R$ action on $\mc{M}_{J_s}(D)$ coming from the $\R$ action on $[0,1]\times\R$, and the \emph{unparametrized moduli space} is denoted by $\wh{\mc{M}_{J_s}}(D)=\mc{M}_{J_s}(D)/\R$. The virtual index bundle of the linearization map $D_u$ gives an element of the $K$-theory of $\mc{B}(D)$. Its dimension is the expected dimension of the moduli space $\mc{M}_{J_s}(D)$, and this dimension is in fact the Maslov index $\mu(D)$, that we had mentioned earlier.  The determinant line bundle of the index bundle, henceforth denoted by $\det(D)$, turns out to be a trivializable line bundle over $\mc{B}(D)$. Therefore, a choice of a nowhere vanishing section on the trivializable line bundle $\det(D)$, produces an orientation of the moduli space $\mc{M}_{J_s}(D)$, and hence an orientation of the unparametrized moduli space $\wh{\mc{M}_{J_s}}(D)$.

\subsection{}\label{para:special}If $D_1\in\mc{D}(x,y)$ and $D_2\in\mc{D}(y,z)$ are domains, then the $2$-chain $D_1+D_2$ lies in $\pi_2(x,z)$. The asymototic behaviors that we had mentioned earlier, along with some globally pre-determined choices, allows us to get a pre-gluing map from $\mc{B}(D_1)\times\mc{B}(D_2)$ to $\mc{B}(D_1+D_2)$. The pullback of the line bundle $\det(D_1+D_2)$ over $\mc{B}(D_1+D_2)$ can be canonically identified with the line bundle $\det(D_1)\wedge \det(D_2)$ over $\mc{B}(D_1)\times\mc{B}(D_2)$ by linearized gluing. An \emph{orientation system} $\mf{o}$ is a choice of a nowhere vanishing section $\mf{o}(D)$ of the line bundle $\det(D)$ for every domain $D\in\mc{D}(x,y)$, and for every pair of generators $x,y\in\mc{G}_{\mc{H}}$, such that if $D_1\in\mc{D}(x,y)$ and $D_2\in\mc{D}(y,z)$, then $\mf{o}(D_1)\wedge\mf{o}(D_2)=\mf{o}(D_1+D_2)$. Therefore, two orientation systems $\mf{o}_1$ and $\mf{o}_2$ disagree on $D_1+D_2$ if and only if they disagree on exactly one of the two domains $D_1$ and $D_2$.

\subsection{}The following describes a method to find all possible orientation systems. Fix a generator $x\in\mc{G}_{\mc{H}}$, and for every other generator $y$, choose a domain $D_y\in\mc{D}(x,y)$. Then choose a set of periodic domains $P_1,\ldots,P_m$, which freely generate $\mc{P}_{\mc{H}}$. Orient the determinant line bundles over the domains $D_y$ and $P_j$ arbitrarily. Since any domain $D\in\mc{D}(y,z)$ can be written uniquely as $D=\sum_j a_j P_j +D_z-D_y$, this choice uniquely specifies an orientation system.
Thus, an orientation system is specified by its values on certain domains $D_y$ and certain periodic domains $P_j$. This allows us to define a chain complex over $\Z$, and it will turn out that the gauge equivalence class of the sign assignment on the chain complex is independent of the orientations of the line bundles $\det(D_y)$. Therefore, declare two orientations systems to be \emph{strongly equivalent} if they agree on all the periodic domains in $\mc{P}_{\mc{H}}$ (or in other words, they agree on all the periodic domains $P_1,\ldots,P_m$). There is a second notion of equivalence, which is of some importance to us, whereby two orientation systems are declared to be \emph{weakly equivalent} if they agree on all the periodic domains in $\mc{P}^0_{\mc{H}}$. Let $\wh{\mc{O}}_{\mc{H}}$ denote the set of weak equivalence classes of orientation systems. Then $\wh{\mc{O}}_{\mc{H}}$ is a torseur over $\Hom(\mc{P}^0_{\mc{H}},\Z/2\Z)$, so there are exactly $2^{b_1(Y)+l-1}$ weak equivalence classes of orientation systems.

\vspace{12 pt}

If $D\in\mc{D}(x,y)$ is a domain, its unparametrized moduli space
$\wh{\mc{M}_{J_s}}(D)$ is a compact, $(\mu(D)-1)$-dimensional manifold
with corners by Gromov compactness and the fact that $J_s$ achieves
transversality; an orientation system $\mf{o}$ determines an
orientation on $\wh{\mc{M}_{J_s}}(D)$.  Therefore, if $\mu(D)=1$, then
$\wh{\mc{M}_{J_s}}(D)$ is a compact oriented zero-dimensional manifold with
corners, or in other words, it is a finite number of signed
points. Let $c(D)$ be the total number of points, counted with
sign. The cornerstone of Floer homology in the present setting, is the
following lemma.

\begin{lem}\cite{POZSz}\label{lem:index2}
If $D\in\mc{D}(x,y)$ is a domain with $\mu(D)=2$, then
$\wh{\mc{M}_{J_s}}(D)$ is an oriented one-dimensional
manifold. Furthermore, if $D=D_1+D_2$, where $D_1\in\mc{D}(x,z)$ and
$D_2\in\mc{D}(z,y)$, with $\mu(D_1)=1$ and $\mu(D_2)=1$, then the
total number of points in the boundary of $\wh{\mc{M}_{J_s}}(D)$ that
correspond to a decomposition of $D$ as $D_1+D_2$, when counted with
signs induced from the orientation of $\wh{\mc{M}_{J_s}}(D)$, equals
$c(D_1)c(D_2)$.
\end{lem}

An immediate corollary is the following: if all the points in the
boundary of $\wh{\mc{M}_{J_s}}(D)$ correspond to such a decomposition --- in other words, if bubbling and boundary degenerations can be ruled out --- then the sum $\sum c(D_1)c(D_2)$ over all such possible decompositions
is zero. This allows us to define the following $(l+1)$-graded chain
complex over $\Z$. This is a well-known chain complex, and it was
first defined by \Ozsvath{} and \Szabo{} for $k=1$. However, for a
general value of $k$, the chain complex was originally not defined
over $\Z$. There are certain subtleties that need to be resolved
before the minus version can be defined over $\Z$, namely, we have to
orient the boundary degenerations in a consistent manner such that the
proofs of Theorems \ref{thm:main}, \ref{thm:secondmain} and
\ref{thm:connectsum} go through; however, those
issues do not appear when we work only in the hat version.

\begin{defn}
  Given an admissible Heegaard diagram $\mc{H}$ for $L$ and an   orientation system $\mf{o}\in\wh{\mc{O}}_{\mc{H}}$, let   $\whcfl_{\mc{H}}(L,\Z,\mf{o})$ be the chain complex freely generated   over $\Z$ by the elements of $\mc{G}_{\mc{H}}$, with the   $(l+1)$ gradings given by $M,A_1,\ldots,A_l$, and the boundary map   given by $\del   x=\sum_{y\in\mc{G}_\mc{H}}\sum_{D\in\mc{D}^0(x,y),\mu(D)=1}c(D)y$.
\end{defn}

\begin{lem}
The map $\del$ on $\whcfl_{\mc{H}}(L,\Z,\mf{o})$ reduces the Maslov grading by $1$, keeps all Alexander gradings fixed, and satisfies $\del^2=0$.
\end{lem}

\begin{proof}
  The claims regarding the gradings follow directly from the
  definitions. To prove that $\del^2=0$, by Lemma \ref{lem:index2}, we
  only need to show that for any empty Maslov index $2$ domain $D$,
  the boundary points of $\wh{\mc{M}}(D)$ do not correspond to
  bubbling or boundary degenerations. However, the shadow of a bubble
  or a boundary degeneration is a $2$-chain in the Heegaard diagram,
  whose boundary lies entirely within the $\alpha$ circles, or
  entirely within the $\beta$ circles. Any such $2$-chain must have
  non-zero coefficient at some $X$ marking, and therefore by
  positivity of domains, the original domain must also have non-zero
  coefficient at that $X$ marking, and therefore, could not have been
  empty.
\end{proof}

Even though we did not specify in the notations, $\whcfl_{\mc{H}}(L,\Z,\mf{o})$ might also depend on the path of almost
complex structures $J_s$ on $\Sym^{g+k-1}(\Sigma_g)$.  However, the
homology $H_*(\whcfl_{\mc{H}}(L,\Z,\mf{o}))$, as an $(l+1)$-graded
object, depends only on the link $L$, the numbers of $X$ markings,
$m_i$, that lie on the $i^{\text{th}}$ link component for each $i$,
and the weak equivalence class of the orientation system $\mf{o}$.

\begin{thm}\label{thm:main}
  For a fixed Heegaard diagram $\mc{H}$ and a fixed path of almost
  complex structures $J_s$, if $\mf{o}_1$ and $\mf{o}_2$ are weakly
  equivalent, then the two chain complexes
  $\whcfl_{\mc{H}}(L,\Z,\mf{o}_1)$ and
  $\whcfl_{\mc{H}}(L,\Z,\mf{o}_2)$ are isomorphic. If $\mc{H}_1$ and
  $\mc{H}_2$ are two different Heegaard diagrams for the same link
  $L$, such that in both $\mc{H}_1$ and $\mc{H}_2$, the
  $i^{\text{th}}$ link component $L_i$ is represented by $m_i$ $X$
  markings and $m_i$ $O$ markings, and if $J_{s,1}$ and $J_{s,2}$ are
  two paths of almost complex structures on the two symmetric
  products, then there is a bijection $f$ between
  $\wh{\mc{O}}_{\mc{H}_1}$ and $\wh{\mc{O}}_{\mc{H}_2}$, such that for
  every $\mf{o}\in\wh{\mc{O}}_{\mc{H}_1}$, the
  homology $H_*(\whcfl_{\mc{H}_1}(L,\Z,\mf{o}))$ is isomorphic to
  the homology $H_*(\whcfl_{\mc{H}_2}(L,\Z,f(\mf{o})))$, as
  $(l+1)$-graded groups.
\end{thm}

\begin{proof}
  This is neither a new type of a theorem, nor a new idea of a proof.
  For the first part, let $\mf{o}_1$ and $\mf{o}_2$ be two weakly
  equivalent orientation systems. We are going to define a map
  $t:\mc{G}_{\mc{H}}\rightarrow\{\pm 1\}$ in the following way.  Call
  two generators $x$ and $y$ to be connected if there is an empty
  domain $D\in\mc{D}^0(x,y)$. For each connected component of
  $\mc{G}_{\mc{H}}$, choose a generator $x$ in that connected
  component, and declare $t(x)=1$. For every other generator $y$ in
  that connected component, choose an empty domain
  $D_y\in\mc{D}^0(x,y)$, and declare $t(y)=1$ if $\mf{o}_1(D_y)$
  agrees with $\mf{o}_2(D_y)$, and $t(y)=-1$ otherwise.  Since
  $\mf{o}_1$ and $\mf{o}_2$ agree on all the empty periodic domains,
  $t$ is a well-defined function. Furthermore, for any empty Maslov
  index $1$ domain $D\in\mc{D}^0(x,y)$, the contribution
  $c_{\mf{o}_1}(D)$ coming from $\mf{o}_1$ is related to the
  contribution $c_{\mf{o}_2}(D)$ coming from $\mf{o}_2$ by the
  equation $c_{\mf{o}_1}(D)=t(x)t(y)c_{\mf{o}_2}(D)$. That shows that
  the two chain complexes are isomorphic via the map $x\mapsto t(x)x$.

  For the second part of the theorem, recall the well known fact that
  if two Heegaard diagrams $\mc{H}_1$ and $\mc{H}_2$ represent the
  same link $L$, such that each component of the link has the same
  number of $X$ and $O$ markings in both the Heegaard diagrams, then
  they can be related to one another by a sequence of isotopies,
  handleslides, and stabilizations. This essentially follows from
  \cite[Proposistion 7.1]{POZSz} and \cite[Lemma
    2.4]{CMPOSS}. However, during the isotopies, we do not require the
  $\alpha$ circles to remain transverse to the $\beta$ circles.
  Therefore, we can assume that $\mc{H}_1$ and $\mc{H}_2$ are related
  by one of the following elementary moves: changing the path of
  almost complex structures $J_s$ by an isotopy $J_{s,t}$; a
  stabilization in a neighborhood of a marked point; a sequence of
  isotopies and handleslides of the $\alpha$ circles in the complement of the
  marked points; or a sequence of isotopies and handleslides of the
  $\be$ circles in the complement of the
  marked points.

For the case of a stabilization, or an isotopy of the
path of almost complex structures, there is a natural identification
between $\mc{P}^0_{\mc{H}_1}$ and $\mc{P}^0_{\mc{H}_2}$, and a natural
identification of the determinant line bundles over the corresponding
empty periodic domains. Since a weak equivalence class of an
orientation system is determined by its values on the empty periodic
domains, this produces a natural identification between
$\wh{\mc{O}}_{\mc{H}_1}$ and $\wh{\mc{O}}_{\mc{H}_2}$. The proof that
the two homologies are isomorphic for the corresponding orientation
systems is immediate for the case of a stabilization, and follows from
the usual arguments of \cite{POZSz} for the other cases. We
do not encounter any new problems, since boundary degenerations are
still ruled out by the marked points.

For the remaining cases, namely, the case of isotopies and
handleslides of $\al$ circles or $\be$ circles, the isomorphism is
established by counting holomorphic triangles. Let us assume that the
$\alpha$ circles are changed to the $\gamma$ circles by a sequence of
isotopies and handleslides in the complement of the marked points. Out
of the $2^{g+k-1}$ weak equivalence classes of orientation systems in
the Heegaard diagram $\mc{H}_3=(\Sigma,\gamma,\alpha,z,w)$, there is a
unique one $\mf{o}_3$, for which the homology of $\mc{H}_3$ is
torsion-free. Each empty periodic domain in $\mc{H}_2$ can be written
uniquely as a sum of empty periodic domains in $\mc{H}_1$ and
$\mc{H}_3$. Therefore, we have a natural bijection between
$\wh{\mc{O}}_{\mc{H}_1}$ and $\wh{\mc{O}}_{\mc{H}_2}$: given an
orientation system $\mf{o}\in\wh{\mc{O}}_{\mc{H}_1}$, we can patch it
with $\mf{o}_3$, to get an orientation system
$f(\mf{o})\in\wh{\mc{O}}_{\mc{H}_2}$. The triangle map, evaluated on
the top generator of the homology of $\mc{H}_3$, provides the required
isomorphism between the homology of $\mc{H}_1$ and the homology of
$\mc{H}_2$, for the corresponding orienation systems. The same proof
from \cite{POZSz} goes through without any problems since we do not encounter any boundary degenerations.
\end{proof}

Let $\vec{m}=(m_1,\ldots,m_l)$. The above theorem shows that
$H_*(\whcfl_{\mc{H}}(L,\Z,\mf{o}))$ is an invariant of the link $L$
inside the three-manifold, a choice of a weak equivalence class of an
orientation system $\mf{o}$, and the vector $\vec{m}$. Let us
henceforth denote the homology as $\whhfl_{\vec{m}}(L,\Z,\mf{o})$. We
now investigate the dependence of
$\whhfl_{\vec{m}}(L,\Z,\mf{o})$ on $\vec{m}$.

\begin{thm}\label{thm:secondmain}
  Let $\mc{H}$ be a Heegaard diagram for a link $L$, where the $i\ith$
  component $L_i$ is represented by $m_i$ $X$ markings and $m_i$ $O$
  markings, and let $\mc{H}'$ be a Heegaard diagram for the same link,
  where $L_i$ is represented by $m'_i=(m_i+\delta_{i_0 i})$ $X$
  markings and $m'_i$ $O$ markings, for some fixed $i_0$.  Then
  there is a bijection $f$ between $\wh{\mc{O}}_{\mc{H}}$ and
  $\wh{\mc{O}}_{\mc{H}'}$ such that for every weak equivalence class
  of orientation system $\mf{o}$, $\whhfl_{\vec{m}'}(L,\Z,\mf{o})$ is
  isomorphic to $\whhfl_{\vec{m}}(L,\Z,f(\mf{o}))\otimes Q_{i_0}$ as
  $(l+1)$-graded groups.
\end{thm}

\begin{proof}
  Consider the Riemann sphere $S$ with one $\alpha$ circle and one
  $\beta$ circle, intersecting each other at two points $p$ and $q$.
  Put two $X$ markings, one $O$ marking and one $W$ marking, one in
  each of the four elementary domains of
  $S\sm(\alpha\cup\beta)$, such that the boundary of either of
  the two elementary domains that contain an $X$ marking runs from
  $p$ to $q$ along the $\alpha$ circle, and from $q$ to $p$ along the
  $\beta$ circle. Remove a small disk in the neighborhood of the point
  $W$. In the Heegaard diagram $\mc{H}$, choose an $X$ marking that
  lies in $L_{i_0}$, and remove a small disk in the neighborhood of
  that point. Then connect the diagram $\mc{H}$ to the sphere $S$ via
  the `neck' $S^1\times [0,T]$ to get a new Heegaard diagram for the
  same link, where $L_i$ is represented by $m'_i$ $X$
  markings, and $m'_i$ $O$ markings. This process is
  shown in Figure \ref{fig:stabilization}.  By Theorem \ref{thm:main},
  we can assume that the new Heegaard diagram is $\mc{H}'$. There is a
  natural correspondance between $\mc{P}^0_{\mc{H}}$ and
  $\mc{P}^0_{\mc{H}'}$, and this induces the bijection $f$ between
  $\wh{\mc{O}}_{\mc{H}}$, and $\wh{\mc{O}}_{\mc{H}'}$.

\begin{figure}
\psfrag{z}{$X$}
\psfrag{w}{$O$}
\psfrag{p}{$p$}
\psfrag{q}{$q$}
\psfrag{a}{$\alpha$}
\psfrag{b}{$\beta$}
\begin{center}
\includegraphics[width=\textwidth]{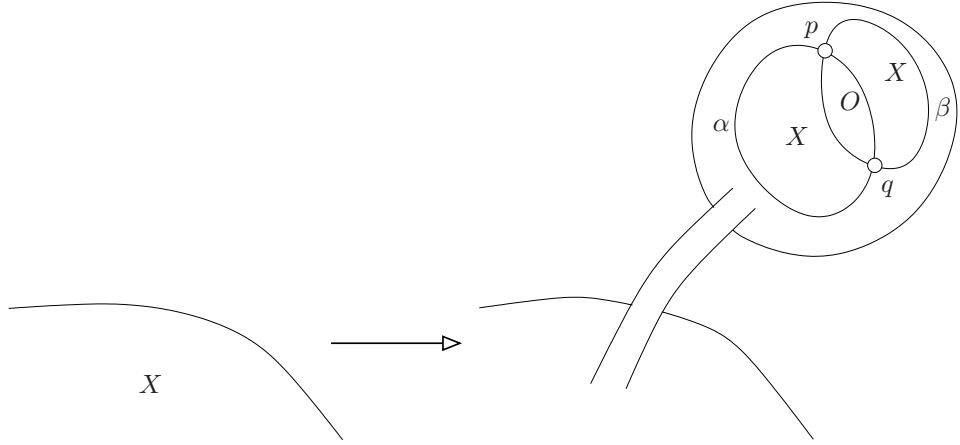}
\end{center}
\caption{The Heegaard diagrams $\mc{H}$ and $\mc{H}'$.}\label{fig:stabilization}
\end{figure}

Fix $\mf{o}\in\wh{\mc{O}}_{\mc{H}}$. As $(l+1)$-graded groups,
$\whcfl_{\mc{H}'}(L,\Z,\mf{o})=\whcfl_{\mc{H}}(L,\Z,f(\mf{o}))\otimes(\Z\oplus\Z)$,
where one $\Z$ corresponds to all the generators that contain the
point $p$, and has $(M,A_1,\ldots,A_l)$ multi-grading
$(0,0,\ldots,0)$, and the other $\Z$ corresponds to all the generators
that contain the point $q$, and has $(M,A_1,\ldots,A_l)$ multi-grading
$(-1,-\delta_{i_0 1},\ldots,-\delta_{i_0 l})$. We simply need to show
that the same identity holds as chain complexes. For this, it is
enough to show that there are no boundary maps from the generators
that contain the point $p$ to the generators that contain the point
$q$.

Following the arguments from \cite{POZSzlinkinvariants}, we extend the
`neck length' $T$, and move the point $W$ close to the $\alpha$ circle
in $S$. After choosing $T$ sufficiently large and $W$ sufficiently
close to the $\alpha$ circle, if there is an empty positive Maslov
index $1$ domain $D$, joining a generator containing $p$ to a
generator containing $q$, such that $c(D)\neq 0$, then $D$ must
correspond to a positive, Maslov index $2$ domain in $\mc{H}$ that
avoids all the $O$ markings and whose boundary lies entirely on the
$\alpha$ circles. However, any non-trivial domain in $\mc{H}$ whose
boundary lies entirely on the $\alpha$ circles must have non-zero
coefficients at some $O$ marking, thus producing a contradiction, and
thereby finishing the proof.
\end{proof}

Henceforth, denote $\whhfl_{(1,\ldots,1)}(L,\Z,\mf{o})$ by $\whhfl(L,\Z,\mf{o})$. 
Theorems \ref{thm:main} and \ref{thm:secondmain} imply:

\begin{thm}\label{thm:fourthmain}
  Let $\mc{H}$ be a Heegaard diagram for a link $L\subset S^3$ with
  $l$ components, such that the $i\ith$ component $L_i$ is represented
  by exactly $m_i$ $X$ markings, and exactly $m_i$ $O$ markings. Then
  the $2^{l-1}$ homology groups $\whhfl_{\vec{m}}(L,\Z,\mf{o})$ are isomorphic to
  the $2^{l-1}$ groups
  $\whhfl(L,\Z,\mf{o})\otimes_i(\otimes^{m_i-1}Q_i)$.
\end{thm}

We are almost done with the construction that we had set out to do. Given a link $L\subset S^3$ with $l$ components, we have produced $2^{l-1}$ $\Z$-valued $(l+1)$-graded homology groups $\whhfl(L,\Z,\mf{o})$. We would like to finish this section by showing that when we combine the $l$ Alexander gradings into one, then we get the $2^{l-1}$ $\Z$-valued bi-graded homology groups $\whhfk(L,\Z,\mf{o})$. Recall that the groups $\whhfk(L,\Z,\mf{o})$ are constructed by viewing the link $L\subset Y$ as a knot in $Y\#^{l-1}(S^1\times S^2)$, and then looking at the knot Floer homology. Therefore, the following lemma is all that we need.

\begin{thm}\label{thm:connectsum}
  Let $\mc{H}$ be a Heegaard diagram for a link $L\subset Y$ with
  $(l+1)$ components, such that each component is represented by one
  $X$ and one $O$ marking. Let $\wt{L}$ be the link with $l$
  components in $Y\#(S^1\times S^2)$, whose $l\ith$ component
  $\wt{L}_l$ is obtained by connect summing $L_{l+1}$ and $L_l$
  through the one-handle, and let $\wt{\mc{H}}$ be a Heegaard diagram
  for $\wt{L}$, where $\wt{L}_i$ is represented by $(1+\delta_{il})$
  $X$ markings and $(1+\delta_{il})$ $O$ markings. Then, there is a
  bijection $f$ between $\wh{\mc{O}}_{\mc{H}}$ and
  $\wh{\mc{O}}_{\wt{\mc{H}}}$, such that for all
  $\mf{o}\in\wh{\mc{O}}_{\mc{H}}$,
  $H_*(\whcfl_{\wt{\mc{H}}}(\wt{L},\Z,f(\mf{o})))=\whhfl(L,\Z,\mf{o})\otimes
  Q_l$ as $(l+1)$-graded groups, where the $(l+1)$ gradings on the
  left hand side are $(M,A_1,\ldots,A_{l-1},A_l+A_{l+1})$.
\end{thm}

\begin{proof}
  This proof is very similar to the proof of Theorem
  \ref{thm:secondmain}. Once more, consider the Riemann sphere $S$
  with one $\alpha$ circle and one $\beta$ circle, intersecting each
  other at two points $p$ and $q$. Put two $X$ markings and two $W$
  marking, one in each of the four elementary domains of
  $S\sm(\alpha\cup\beta)$, such that the boundary of either of
  the two elementary domains that contain an $X$ marking runs from
  $p$ to $q$ along the $\alpha$ circle, and from $q$ to $p$ along the
  $\beta$ circle. Remove two small disks in the neighborhoods of the
  $W$ markings. In the Heegaard diagram $\mc{H}$, remove two small
  disks in the neighborhoods of the the two $X$ markings
  that lie in $L_l$ and $L_{l+1}$. Then connect $\mc{H}$ to the sphere
  $S$ via the two `necks,' $S^1\times[0,T_1]$ and $S^1\times[0,T_2]$,
  as shown in Figure \ref{fig:connectsum}. The resulting picture is a
  Heegaard diagram for the link $\wt{L}\subset Y\#(S^1\times S^2)$,
  where the $i\ith$ component $\wt{L}_i$ is represented by
  $(1+\delta_{il})$ $X$ markings and $(1+\delta_{il})$ $O$ markings.
  By the virtue of Theorem \ref{thm:main}, we can assume that this
  Heegaard diagram is $\wt{\mc{H}}$.

\begin{figure}
\psfrag{z}{$X$}
\psfrag{p}{$p$}
\psfrag{q}{$q$}
\psfrag{a}{$\alpha$}
\psfrag{b}{$\beta$}
\begin{center}
\includegraphics[width=\textwidth]{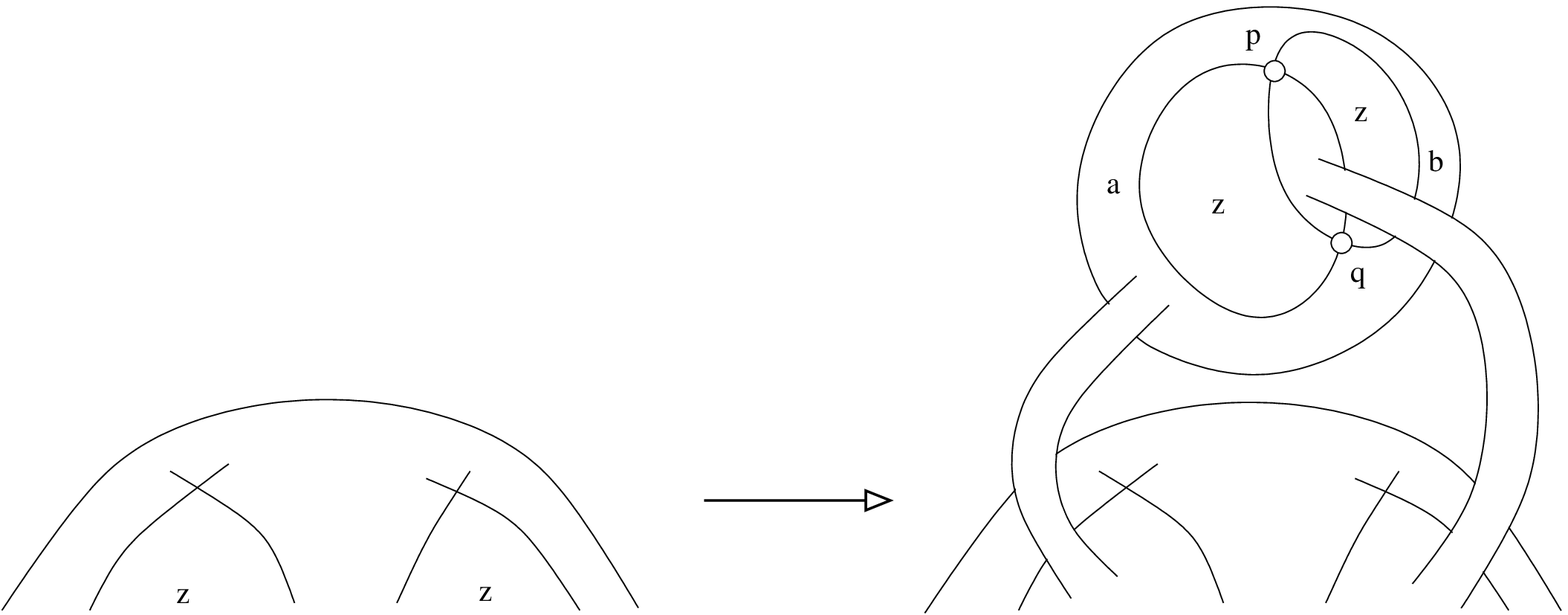}
\end{center}
\caption{The Heegaard diagrams $\mc{H}$ and $\wt{\mc{H}}$.}\label{fig:connectsum}
\end{figure}

An empty periodic domain in $\mc{H}$ gives rise to an empty periodic
domain in $\wt{\mc{H}}$. In the other direction, an empty periodic
domain in $\wt{\mc{H}}$ gives rise to a periodic domain in $\mc{H}$
which does not pass through any of the $O$ markings. Since each
component of the link $L$ is null-homologous in $Y$, such a periodic
domain is an empty periodic domain. Therefore, there is a natural
correspondance between the empty periodic domains of $\mc{H}$ and
$\wt{\mc{H}}$, and this induces the bijection $f$ between
$\wh{\mc{O}}_{\mc{H}}$ and $\wh{\mc{O}}_{\wt{\mc{H}}}$.

Fix $\mf{o}\in\wh{\mc{O}}_{\mc{H}}$. It
is immediate that as $(l+1)$-graded groups,
$\whcfl_{\wt{\mc{H}}}(\wt{L},\Z,f(\mf{o}))=\whcfl_{\mc{H}}(L,\Z,\mf{o})\otimes Q_l$.
However, quite like the case of Theorem \ref{thm:secondmain}, for
sufficiently large `neck lengths' $T_1$ and $T_2$, and with the two
$W$ markings sufficiently close to the $\alpha$ circle on $S$, the
above identity holds even as chain complexes.
\end{proof}

Before we conclude this section, a note regarding absolute gradings is
due. So far, we have worked with relative Maslov grading and relative
Alexander gradings.  However, for links in $S^3$, and for links in
$\#^m (S^1\times S^2)$ that we obtain from links in $S^3$ by the
connect sum process described in Theorem \ref{thm:connectsum}, there
is a well-defined way to lift these gradings to absolute gradings, as
defined in \cite[Theorem 7.1]{POZSz4manifolds}, \cite[Subsections 3.3
  and 3.4]{POZSzknotinvariants} and \cite[Lemma 4.6 and Equation
  24]{POZSzlinkinvariants}. Since this is an oft-studied scenario, for
such links, let us improve the earlier theorems, and
henceforth work with absolute gradings.

\begin{lem}\label{lem:absolutegrading}
  For links in $\#^m(S^1\times S^2)$ that come from links in $S^3$ by
  the connect sum operation as decribed in Theorem
  \ref{thm:connectsum}, the isomorphisms in Theorems \ref{thm:main},
  \ref{thm:secondmain}, \ref{thm:fourthmain} and \ref{thm:connectsum}
  preserve the absolute gradings.
\end{lem}

\begin{proof}
  Recall that the isomorphisms in question come from chain maps that
  preserve the relative gradings. Therefore, each such chain map must
  shift each absolute grading by a fixed integer on the entire
  chain complex. We want to show that each of these shifts is zero.

  Since the absolute gradings are defined on the generators
  themselves, this shift is unchanged if instead of working over $\Z$,
  we tensor everything with $\F_2$ and work over $\F_2$. However,
  since the Heegaard Floer homology of $\#^m(S^1\times S^2)$ is
  non-trivial over $\F_2$, in each case, the homology of the entire
  chain complex is non-trivial over $\F_2$. Furthermore, the 
  maps induced on the homology over $\F_2$ preserve the absolute gradings
  \cite{POZSz4manifolds}, \cite{POZSzknotinvariants},
    \cite{POZSzlinkinvariants}. Therefore, all the shifts are zero, and each
  of the chain maps preserves all the gradings.
\end{proof}

\section{Grid diagrams}\label{sec:griddiagrams}

A \emph{planar grid diagram of index $N$} is the square
$S=[0,N]\times[0,N]\subset\R^2$, with the following additional
structures: if $1\leq i\leq N$, the horizontal line $y=(i-1)$ is
called $\alpha_i$, the $i\ith$ $\alpha$ arc, and the vertical line
$x=(i-1)$ is called $\beta_i$, the $i\ith$ $\beta$ arc; there are $2N$
markings, denoted by $X_1,\ldots,X_N,O_1,\ldots,O_N$, such that each
component of $S\sm(\bigcup_i\alpha_i)$ contains one $X$ marking and
one $O$ marking, and each component of $S\sm(\bigcup_i\beta_i)$
contains one $X$ marking and one $O$ marking.

A \emph{toroidal grid diagram of index $N$} is obtained from a planar
grid diagram of the same index by identifying the opposite sides of
the square $S$ to form a torus $T$. A careful reader will immediately
observe that this creates a Heegaard diagram $\mc{H}$ for some link
$L$ in $S^3$, and for the rest of the section, we will work with this
Heegaard diagram. The $\alpha$ arcs and the $\beta$ arcs become full
circles, and they are the $\alpha$ circles and the $\beta$ circles
respectively; the $N$ components of $T\sm(\bigcup_i\alpha_i)$ are called
the \emph{horizontal annuli}, and each of them contains one $X$
marking and one $O$ marking; the horizontal annulus with $\alpha_i$ as
the circle on the bottom is called the $i\ith$ horizontal annulus, and
is denoted by $H_i$; the $N$ components of $T\sm(\bigcup_i\beta_i)$ are
called the \emph{vertical annuli}, and each of them also contains one
$X$ marking and one $O$ marking; the vertical annulus with $\beta_i$
as the circle on the left is called the $i\ith$ vertical annulus,
and is denoted by $V_i$; the $N^2$ components of
$T\sm\bigcup_i(\alpha_i\cup\beta_i)$ are the elementary
domains. Therefore, the link $L$ that the toroidal grid diagram
represents, can be obtained in the following way. We assume that the
toroidal grid diagram comes from a planar grid diagram on the square
$S$. Then in each component of $S\sm(\bigcup_i\alpha_i)$, we join the $X$
marking to the $O$ marking by an embedded arc, and in each component
of $S\sm(\bigcup_i\beta_i)$, we join the $O$ marking to the $X$ marking
by an embedded arc, and at every crossing, we declare the arc that
joins $O$ to $X$ to be the overpass. Henceforth, we also assume that
the link $L$ has $l$ components, and the $i\ith$ component $L_i$ is
represented by $m_i$ $X$ markings and $m_i$ $O$ markings, and $\sum_i
m_i =N$.

There is only one $Spin^C$ structure, so generators in $\mc{G}_{\mc{H}}$
correspond to the permuatations in $\mf{S}_N$ as follows: a generator
$x=(x_1,\ldots,x_N)\in\mc{G}_{\mc{H}}$ comes from the permutation
$\sigma\in\mf{S}_n$, where $x_i=\alpha_i\cap\beta_{\sigma(i)}$ for
each $1\leq i\leq N$. The $N$ points $x_1,\ldots,x_N$ are called the
\emph{coordinates} of the generator $x$.

Let $\mf{j}$ be the complex structure on $T$ induced from the standard
complex structure on $S\subset\C$, and let $J_s$ be the constant path
of almost complex structure $\Sym^N(\mf{j})$ on $\Sym^N(T)$. After a
slight perturbation of the $\al$ and the $\be$ circles, we can ensure
that $J_s$ achieves transversality for all domains up to Maslov index
two \cite[Lemma 3.10]{RL}. Henceforth, we work with these perturbed
$\al$ and $\be$ circles and this path of nearly symmetric almost
complex structure.

Consider the $2^{l-1}$ chain complexes $\whcfl_{\mc{H}}(L,\Z,\mf{o})$. The boundary maps in each of the chain complexes correspond to objects called \emph{rectangles}. A rectangle $R$ joining a generator $x$ to a generator $y$ is a $2$-chain generated by the elementary domains of $\mc{H}$, such that the following conditions are satisfied: $R$ only has coefficients $0$ and $1$; the closure of the union of the elementary domains where $R$ has coefficient $1$ is a disk embedded in $T$ with four corners, or in other words, it looks like a rectangle; the top-right corner and the bottom-left corner of $R$ are coordinates of $x$; the top-left corner and the bottom-right corner of $R$ are coordinates of $y$; the generators $x$ and $y$ share $(N-2)$ coordinates; and $R$ does not contain any coordinates of $x$ or any coordinates of $y$ in its interior. It is easy to check that the rectangles are precisely the positive Maslov index one domains. We denote the set of all rectangles joining $x$ to $y$ by $\mc{R}(x,y)\subset \mc{D}(x,y)$. The set $\mc{R}(x,y)$ is empty unless $x$ and $y$ differ in exactly two coordinates, and even then, $\left|\mc{R}(x,y)\right|\leq 2$.

\begin{lem}\cite[Theorem 1.1]{CMPOSS}\label{lem:maslovone}
If $D\in\mc{D}(x,y)$ is a domain with $\mu(D)\leq 0$, then the
unparametrized moduli space $\wh{\mc{M}_{J_s}}(D)$ is empty. If
$D\in\mc{D}(x,y)$ is a Maslov index one domain such that
$\wh{\mc{M}_{J_s}}(D)$ is non-empty, then $D$ is a
rectangle. Conversely, if $R\in\mc{R}(x,y)$ is a rectangle, then
$\wh{\mc{M}_{J_s}}(R)$ consists of exactly one point, and hence
$\left|c(R)\right|=1$.
\end{lem}

If $D\in\mc{D}(x,y)$, we say that $D$ can be decomposed as a sum of
two rectangles if there exists a generator $z\in\mc{G}_{\mc{H}}$ and
rectangles $R_1\in\mc{R}(x,z)$ and $R_2\in\mc{R}(z,y)$ such that
$D=R_1+R_2$. It is easy to check that the domains that can
be decomposed as sum of two rectangles are precisely the
positive Maslov index two domains. For any generator $x\in\mc{G}_T$,
there are exactly $2N$ Maslov index two positive domains in
$\mc{D}(x,x)$, namely the ones coming from the horizontal annuli
$H_1,\ldots,H_N$ and the vertical annuli $V_1,\ldots,V_N$.

\begin{lem}\label{lem:maslovtwo}
  If $D\in\mc{D}(x,y)$ is a Maslov index two domain such that
  $\wh{\mc{M}_{J_s}}(D)$ is non-empty, then $D$ can be decomposed as a
  sum of two rectangles. Conversely, if $D\in\mc{D}(x,y)$ can be
  decomposed as a sum of two rectangles, then $\wh{\mc{M}_{J_s}}(D)$
  is a compact $1$-dimensional manifold with exactly two
  endpoints. Furthermore, if $x=y$ (i.e. if $D$ comes from a
  horizontal or a vertical annulus), then one of the endpoints
  corresponds to the unique way of decomposing $D$ as a sum of two
  rectangles, while the other endpoint corresponds to an $\alpha$ or a
  $\beta$ boundary degeneration; and if $x\neq y$, then $D$ can be
  decomposed as a sum of two rectangles in exactly two ways, and the
  two endpoints correspond to the two decompositions.
\end{lem}

Lemma \ref{lem:maslovone} implies that once we choose an orientation system $\mf{o}$ (and not just a weak equivalence class of orientation systems), we get a function $c_{\mf{o}}$ from the set of all rectangles to $\{-1,1\}$. Lemma \ref{lem:maslovtwo} in conjunction with Lemma \ref{lem:index2} implies that if a domain $D\in\mc{D}(x,y)$ can be decomposed as a sum of two rectangles in two different ways $D=R_1+R_2=R_3+R_4$, then $c_{\mf{o}}(R_1)c_{\mf{o}}(R_2)=-c_{\mf{o}}(R_3)c_{\mf{o}}(R_4)$. This naturally leads to the definition of a sign assignment.

\begin{defn}
  A \emph{sign assignment} $s$ is a function from the set of all rectangles   to the set $\{-1,1\}$, such that the following condition is   satisfied: if $x,y,z,z'\in\mc{G}_{\mc{H}}$ are distinct generators,   and if $R_1\in\mc{R}(x,z)$, $R_2\in\mc{R}(z,y)$,   $R'_1\in\mc{R}(x,z')$, $R'_2\in\mc{R}(z',y)$ are rectangles with   $R_1+R_2=R'_1+R'_2$, then $s(R_1)s(R_2)=-s(R'_1)s(R'_2)$. Two sign   assignments $s_1$ and $s_2$ are said to be \emph{gauge equivalent} if there   is a function $t:\mc{G}_{\mc{H}}\rightarrow\{-1,1\}$, such that   $s_1(R)=t(x)t(y)s_2(R)$, for all $x,y\in\mc{G}_{\mc{H}}$ and for all   $R\in\mc{R}(x,y)$.
\end{defn}

In particular, a true sign assignment, as defined in \cite[Definition 4.1]{CMPOZSzDT}, is a sign assignment. Let $f$ be the map from the set of all orientation systems to the set of all sign assignments such that for all rectangles $R$, $f(\mf{o})(R)=c_{\mf{o}}(R)$. In this section, we will show that there are exactly $2^{2N-1}$ gauge equivalence classes of sign assignments on the grid diagram. We will put a weak equivalence on the sign assignments, which is weaker than the gauge equivalence. We will prove that there are exactly $2^{l-1}$ weak equivalence classes of sign assignments, and the map $f$ induces a bijection $\wt{f}$ between the set of weak equivalence classes of orientation systems and the set of weak equivalence classes of sign assignments. This will allow us to combinatorially calculate $\whcfl_{\mc{H}}(L,\Z,\mf{o})$ for all $\mf{o}\in\wh{\mc{O}}_{\mc{H}}$, and thereby calculate $\whhfl(L,\Z)$ in all the $2^{l-1}$ versions. As a corollary, this will also show that any sign assignment (in particular, the one constructed in \cite{CMPOZSzDT}) computes $\whhfl(L,\Z,\mf{o})$ for some orientation system $\mf{o}$.

We have an explicit (although slightly artificial) correspondance between the generators in $\mc{G}_{\mc{H}}$ and the elements of the symmetric group $\mf{S}_N$, whereby a permutation $\sigma\in\mf{S}_N$ gives rise to the generator $x=(x_1,\ldots,x_N)$ with $x_i=\alpha_i\cap\beta_{\sigma(i)}$. There is the following very natural partial order on the permutations: a reduction of a permutation $\tau$ is a permuation obtained by pre-composing $\tau$ by some transposition $(i,j)$ where $i<j$ and $\tau(i)>\tau(j)$; the permutation $\sigma$ is declared to be smaller than the permutation $\tau$, if $\sigma$ can be obtained from $\tau$ by a sequence of reductions. This induces a partial order $\prec$ on the elements of $\mc{G}_{\mc{H}}$.

For $x,y\in\mc{G}_{\mc{H}}$, if $y\prec x$ and there does not exist any $z\in\mc{G}_{\mc{H}}$ such that $y\prec z\prec x$, then we say that $x$ \emph{covers} $y$, and write that as $y\leftarrow x$. If we view the toroidal grid diagram as one coming from a planar grid diagram on $S=[0,N]\times[0,N]$, then $y\leftarrow x$ precisely when there is a rectangle from $x$ to $y$ contained in the subsquare $S'=[0,N-1]\times[0,N-1]$.

The poset $(\mc{G}_{\mc{H}},\prec)$ is a well-understood object \cite{shellPHE}. There is a unique minimum $p\in\mc{G}_{\mc{H}}$, which corresponds to the identity permutation. In particular, the Hasse diagram of $(\mc{G}_{\mc{H}},\prec)$, viewed as an unoriented graph, is connected. There is a unique maximum $q\in\mc{G}_{\mc{H}}$, which corresponds to the permutation that maps $i$ to $(N+1-i)$. The poset is shellable, which means that there is a total ordering $<$ on the maximal chains, such that if $\mf{m}_1$ and $\mf{m}_2$ are two maximal chains with $\mf{m}_1<\mf{m}_2$, then there exists a maximal chain $\mf{m}_3<\mf{m}_2$ with $\mf{m}_1\cap\mf{m}_2\subseteq\mf{m}_3\cap\mf{m}_2=\mf{m}_2\sm\{z\}$ for some $z\in\mf{m}_2$. This in particular implies that given any two maximal chains $\mf{m}_1$ and $\mf{m}_2$, we can get from $\mf{m}_2$ to $\mf{m}_1$ via a sequence of maximal chains, where we get from one maximal chain to the next by changing exactly one element.

Given a sign assignment $s$ and a generator $x\in\mc{G}_{\mc{H}}$, we define two functions $h_{s,x},v_{s,x}:\{1,\ldots,N\}\rightarrow\{-1,1\}$, called the \emph{horizontal function} and the \emph{vertical function}, as follows: let $D\in\mc{D}(x,x)$ be Maslov index two positive domain which corresponds to the horizontal annulus $H_i$; then, $D$ can be decomposed as a sum of two rectangles in a unique way, and define the horizontal function $h_{s,x}(i)$ as the product of the signs of the two rectangles. The vertical function $v_{s,x}(i)$ is constructed similarly by considering the vertical annulus $V_i$ instead. Clearly, the horizontal and the vertical functions depend only on the gauge equivalence class of the sign assignment. The following theorem shows that the functions do not depend on the choice of the generator $x$, and will henceforth be denoted by $h_s$ and $v_s$.

\begin{thm}\label{thm:projection}
For any sign assignment $s$, for any two generators $x,y\in\mc{G}_{\mc{H}}$, and for any $1\leq i\leq N$, the horizontal and the vertical functions satisfy $h_{s,x}(i)=h_{s,y}(i)$ and $v_{s,x}(i)=v_{s,y}(i)$.
\end{thm}

\begin{proof}
Fix a sign assignment $s$, and fix $i\in\{1,\ldots,N\}$. We will only prove the statement for the vertical function; the argument for the horizontal function is very similar. Given $z\in\mc{G}_{\mc{H}}$, let $(z',R_z,R'_z)$ be the unique triple with $z'\in\mc{G}_T$, $R_z\in\mc{R}(z,z')$ and $R'_z\in\mc{R}(z',z)$ such that $R_z+R'_z\in\mc{D}(z,z)$ comes from the vertical annulus $V_i$. We simply want to show that for any two generators $x,y\in\mc{G}_{\mc{H}}$, $s(R_x)s(R'_x)=s(R_y)s(R'_y)$. Recall the partial order on $\mc{G}_{\mc{H}}$. The corresponding Hasse diagram, when viewed as an unoriented graph, is connected; therefore, it is enough to prove the above statement when $y\leftarrow x$. Thus, we can assume that there exists a rectangle $R\in\mc{R}(x,y)$. We end the proof by considering the following two cases.

\begin{figure}
\begin{center}
\includegraphics[height=170pt]{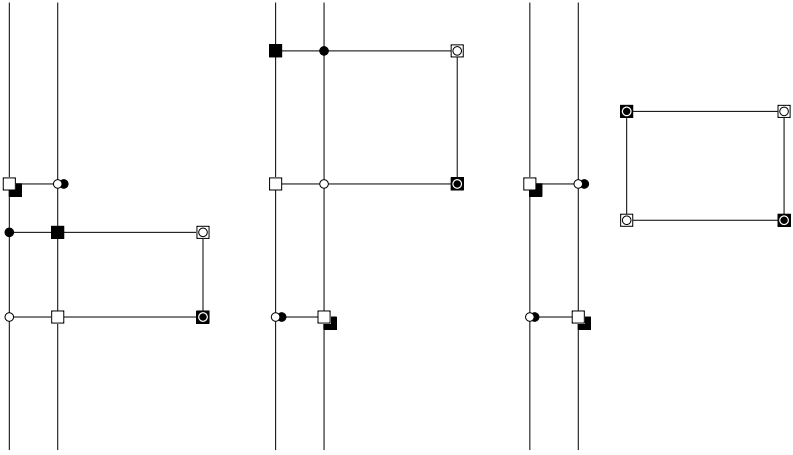}
\end{center}
\caption{The case when $y$ and $x'$ disagree in exactly $3$ or exactly $4$ coordinates. The coordinates of $x$, $y$, $x'$ and $y'$ are denoted by white circles, black circles, white squares and black squares, respectively.}\label{fig:projection}
\end{figure}

\emph{The generators $y$ and $x'$ disagree on none of the   coordinates.}  In this case, $y=x'$, $y'=x$, $R_x=R'_y$ and $R_y=R'_x$.  The equality $s(R_x)s(R'_x)=s(R_y)s(R'_y)$ follows trivially.

\emph{The generators $y$ and $x'$ disagree on exactly three or exactly   four coordinates.} In this case, there exists a rectangle $R'\in\mc{R}(x',y')$, such that $R_x+R'=R+R_y\in\mc{D}(x,y')$ and $R'_x+R=R'+R'_y\in\mc{D}(x',y)$. The three essentially different types of diagrams that might appear (up to a rotation by $180^{\circ}$) are illustrated in Figure \ref{fig:projection}. Therefore,  $s(R_x)s(R')=-s(R)s(R_y)$ and $s(R'_x)s(R)=-s(R')s(R'_y)$. Multiplying, we get the required identity $s(R_x)s(R'_x)=s(R_y)s(R'_y)$.
\end{proof}

The following two theorems will establish that there are exactly $2^{2N-1}$ gauge equivalence classes of sign assignments. Let $\Phi$ be the map from the set of gauge equivalence classes of sign assignments to $\{-1,1\}^{2N-1}$ given by $s\rightarrow (h_s(1),\ldots,\allowbreak h_s(N),\allowbreak v_s(1),\ldots,v_s(N-1))$.

\begin{thm}\label{thm:equivalence}
Given functions $g_h,g_v:\{1,\ldots,N\}\rightarrow\{-1,1\}$, such that $\left|g^{-1}_v(1)\right|\equiv \left|g^{-1}_h(-1)\right|\pmod{2}$, there exists a sign assignment $s$, such that $g_h=h_s$ and $g_v=v_s$. Therefore, in particular, the function $\Phi$ from the set of gauge equivalence classes of sign assignments to $\{-1,1\}^{2N-1}$ is surjective.
\end{thm}

\begin{proof}
  By \cite[Theorem 4.2]{CMPOZSzDT}, there exists a sign assignment   $s_0$ such that $h_{s_0}(i)=1$ and $v_{s_0}(i)=-1$ for all   $i\in\{1,\ldots,N\}$. Given   $g_h,g_v:\{1,\ldots,N\}\rightarrow\{-1,1\}$ with   $\left|g^{-1}_v(1)\right|\equiv \left|g^{-1}_h(-1)\right|\pmod{2}$,   we would like to modify $s_0$ to get $s$, such that $g_h=h_s$ and   $g_v=v_s$.

  The general method that we employ to modify a sign assignment $s_1$   to get another sign assignment $s_2$, is the following: we start   with a multiplicative $2$-cochain $m$ which assigns elements of $\{-1,1\}$ to the   elementary domains; if $D$ is a $2$-chain generated by the   elementary domains, then $\langle m,D\rangle$ is simply the   evaluation of $m$ on $D$; then, for a rectangle   $R\in\mc{R}(x,y)$, we define $s_2(R)$ to be $s_1(R)\langle   m,R\rangle$. It is easy to see that $s_2$ is a sign assignment if   and only if $s_1$ is a sign assignment.

  We prove the statement by an induction on the number   $n(g_v,g_h)=\frac{1}{2}(\left|g^{-1}_v(1)\right|+\left|g^{-1}_h(-1)\right|)$.   For the base case, when $n(g_v,g_h)=0$, we can simply choose   $s=s_0$.

  Assuming that the induction hypothesis is proved for $n=k$, let   $g_h,g_v:\{1,\ldots,N\}\allowbreak\rightarrow\{-1,1\}$ be functions with   $n(g_v,g_h)=k+1$. Choose functions   $\wt{g}_h,\wt{g}_v:\{1,\ldots,N\}\allowbreak\rightarrow\{-1,1\}$ such that   $n(\wt{g}_v,\wt{g}_h)=k$ and $\left|\{i\mid     g_v(i)\neq\wt{g}_v(i)\}\right|+\left|\{i\mid     g_h(i)\neq\wt{g}_h(i)\}\right|=2$. By induction, there is a sign   assignment $\wt{s}$ such that $\wt{g}_h=h_{\wt{s}}$ and   $\wt{g}_v=v_{\wt{s}}$. If $\left|\{i\mid     g_v(i)\neq\wt{g}_v(i)\}\right|=2$, consider the two vertical   annuli corresponding to the two values where $g_v$ disagrees with   $\wt{g}_v$, choose a horizontal annulus, and let $m$ be the   $2$-cochain which assigns $(-1)$ to the two elementary domains where   the horizontal annulus intersects the two vertical annuli, and $1$   to every other elementary domain. Similarly, if $\left|\{i\mid     g_h(i)\neq\wt{g}_h(i)\}\right|=2$, consider the two horizontal   annuli corresponding to the two values where $g_h$ disagrees with   $\wt{g}_h$, choose a vertical annulus, and let $m$ be the   $2$-cochain which assigns $(-1)$ to the two elementary domains where   the vertical annulus intersects the two horizontal annuli, and $1$   to every other elementary domain. Finally, if $\left|\{i\mid     g_v\neq\wt{g}_v(i)\}\right|=\left|\{i\mid     g_h\neq\wt{g}_h(i)\}\right|=1$, consider the vertical annulus   corresponding to the value where $g_v$ disagrees with $\wt{g}_v$,   consider the horizontal annulus corresponding to the value where   $g_h$ disagrees with $\wt{g}_h$, and let $m$ be the $2$-cochain   which assigns $(-1)$ to the elementary domain where the vertical   annulus intersects the horizontal annulus, and $1$ to every other   elementary domain. Let $s$ be the sign assignment obtained from   $\wt{s}$ by modifying it by the $2$-cochain $m$. It is fairly   straightforward to check that $g_h=h_s$ and $g_v=v_s$.
\end{proof}

\begin{thm}\label{thm:uniqueness}
The function $\Phi$ from the set of gauge equivalence classes of sign assignments to $\{-1,1\}^{2N-1}$ is injective.
\end{thm}

\begin{proof}
For this proof, we will closely follow the corresponding proof from \cite{CMPOZSzDT}. However, that proof uses the permutahedron whose $1$-skeleton is the Cayley graph of the the symmetric group, where the generators are the adjacent transpositions. In our proof, we will use a different simplicial complex, which is the order complex of the partial order $\prec$ on $\mc{G}_{\mc{H}}$.

Recall that the poset has a unique minimum $p$, and a unique maximum $q$. View the Hasse diagram of the poset as an oriented graph $\mf{g}$. Choose a maximal tree $\mf{t}$ with $p$ as a root, i.e. given any vertex $x$, there is a (unique) oriented path from $p$ to $x$ in $\mf{t}$. The edges of $\mf{g}$ correspond to the rectangles that are supported in $[0,N-1]\times[0,N-1]$. A sign assignment endows the edges of $\mf{g}$ with signs $\pm 1$.

Let us choose a $(2N-1)$-tuple in $\{-1,1\}^{2N-1}$, and let $s$ be a   sign assignment such that the $(2N-1)$-tuple equals $\Phi(s)$. We would like to show that the gauge equivalence class of $s$ is determined. Since $\mf{t}$ is a tree, by replacing the sign assignment $s$ by a gauge equivalent one if necessary, we can assume that $s$ labels all the edges of $\mf{t}$ with $1$'s. We will show that the values of $s$ on all the other edges are now determined.

Now consider any other edge $y\leftarrow x$ in $\mf{g}$. Let $\mf{c}_1$ be the unique oriented path from $p$ to $x$ in $\mf{t}$, and let $\mf{c}_2$ be the unique oriented path from $p$ to $y$ in $\mf{t}$. Choose an oriented path $\mf{c}_0$ from $x$ to $q$ in $\mf{g}$. Let $\mf{m}_1$ be the union of $\mf{c}_1$ and $\mf{c}_0$, and let $\mf{m}_2$ be the union of $\mf{c}_2$, the edge from $y$ to $x$, and $\mf{c}_0$; these can be seen as maximal chains in $(\mc{G}_{\mc{H}},\prec)$. Clearly, $($the product of the signs on the edges in $\mf{m}_1)\cdot($the product of the signs on the edges in $\mf{m}_2)=($the product of the signs on the edges in $\mf{c}_1)\cdot($the product of the signs on the edges in $\mf{c}_2)\cdot($the sign on the edge from $y$ to $x)$. Since $\mf{c}_1\cup\mf{c}_2\subseteq\mf{t}$, the signs on the edges of $\mf{c}_1$ and $\mf{c}_2$ are all $1$, so the sign on the edge from $y$ to $x$ equals $($the product of the signs on the edges in $\mf{m}_1)\cdot($the product of the signs on the edges in $\mf{m}_2)$. Since $(\mc{G}_{\mc{H}},\prec)$ is shellable, $\mf{m}_2$ can be turned into $\mf{m}_1$ through maximal chains by modifying one element at a time. Changing exactly one element of exactly one of the maximal chains negates the above product, so the product depends only on the graph $\mf{g}$. Thus, $s$ is determined on all the edges of $\mf{g}$.

Therefore, we have shown that there exists at most one sign assignment, up to gauge equivalence, on the rectangles that lie in the subsquare $S'=[0,N-1]\times[0,N-1]$. In fact, shellability of our poset also implies that there exists a sign assignment, but we do not need it. The rest of the proof for uniqueness is very similar to the proof from \cite{CMPOZSzDT}, but for the reader's convenience, we repeat the argument. Let $S''\subset T$ be the annular subspace corresponding to the rectangle $[0,N-1]\times[0,N]$ in the planar grid diagram. Next, we show that the value of $s$ is determined on all the rectangles that lie in $S''$.  

This is done by an induction on the (horizontal) width of the rectangles. For the base case, if $R\in\mc{R}(x,y)$ is a rectangle of width one which is not supported in $S'$, then let $R'\in\mc{R}(y,x)$ be the unique rectangle such that $R+R'$ is a vertical annulus. The vertical function $v_s$ determines the product of the signs $s(R)s(R')$, and thereby the sign $s(R)$.

\begin{figure}
\begin{center}
\includegraphics[height=170pt]{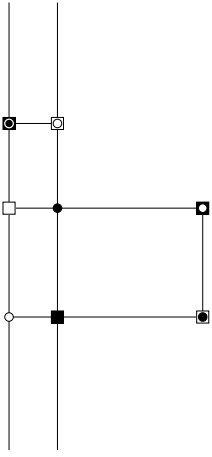}
\end{center}
\caption{The induction step. The coordinates of $x$, $y$,
$y'$ and $z$ are denoted by white circles, white squares, black
squares and black circles, respectively.}\label{fig:uniqueness}
\end{figure}

Assuming that we have proved the uniqueness of sign assignments for all
the rectangles up to width $k$, let $R\in\mc{R}(x,y)$ be a width
$(k+1)$ rectangle. Let $R_1\in\mc{R}(y,z)$ be the width one rectangle
such that the bottom-left corner of $R_1$ is the top-left corner of
$R$. Then there exists a generator $y'\neq y$, a width one rectangle
$R'\in\mc{R}(x,y')$ and a width $k$ rectangle $R'_1\in\mc{R}(y',z)$,
such that $R+R_1=R'+R'_1\in\mc{D}(x,z)$. The situation is illustrated
in Figure \ref{fig:uniqueness}. By induction, the value of
$s$ is determined on $R_1$, $R'$ and $R'_1$. However,
$s(R)s(R_1)=-s(R')s(R'_1)$, and this determines the sign $s(R)$. This
completes the induction and shows that the value of the sign
assignment $s$ is fixed on all the rectangles that are supported
in $S''$. A similar argument, but with the diagrams rotated by
$90^{\circ}$, shows that the value of $s$ is, in fact, determined on
all the rectangles. This completes the proof of uniqueness.
\end{proof}

\begin{lem}\label{lem:product}
For any sign assignment $s$, the product $\prod_{i=1}^N h_s(i)v_s(i)$ equals $(-1)^N$.
\end{lem}

\begin{proof}
By Theorem \ref{thm:equivalence}, there exists a sign assignment $s'$ such that $h_{s'}=h_{s}$, $v_{s'}(i)=v_s(i)$ for $i\in\{1,\ldots,N-1\}$ and $v_{s'}(N)= (-1)^N h_s(N)\prod_{i=1}^{N-1} h_s(i)v_s(i)$. Since $\Phi(s)=\Phi(s')$, by Theorem \ref{thm:uniqueness}, $s$ and $s'$ are gauge equivalent. Therefore, $\prod_{i=1}^N h_s(i)v_s(i)=\prod_{i=1}^N h_{s'}(i)v_{s'}(i)=(-1)^N$.
\end{proof}

Fix a sign assignment $s$ and fix a link component $L_i$. Let $V(L_i)=\{j\mid \text{the }X\text{ marking in }V_j\text{ is in }L_i\}$ and let $H(L_i)=\{j\mid \text{the }X\text{ marking in }H_j\text{ is in }L_i\}$.  The product $(\prod_{j\in H(L_i)}h_s(j))(\prod_{j\in V(L_i)}(-v_s(j)))$ is defined to be the \emph{sign of the link component $L_i$} and is denoted by $r_s(L_i)$.

Call two sign assignments $s_1$ and $s_2$ \emph{weakly equivalent} if $r_{s_1}$ agrees with $r_{s_2}$ on each of the link components. Clearly, if two sign assignments are gauge equivalent, then they are weakly equivalent. Due to Lemma \ref{lem:product}, the product of the signs of all the link components is $1$, and this is the only restriction on these numbers $r_s(L_i)$. Therefore, there are exactly $2^{l-1}$ weak equivalence classes of sign assignments. The following observation yields a direct proof that the chain complex $\whcfl_{\mc{H}}(L,\Z,\mf{o})$ depends only on the weak equivalence class of the sign assignment $f(\mf{o})$.

\begin{lem}
If two sign assignments $s_1$ and $s_2$ are weakly equivalent, then there exists a sign assignment $s'_2$, which is gauge equivalent to $s_2$, such that $s_1$ and $s'_2$ agree on all the rectangles that avoid the $X$ markings and the $O$ markings.
\end{lem}

\begin{proof}
Since $s_1$ and $s_2$ are weakly equivalent, a proof similar to the proof of Theorem \ref{thm:equivalence} shows that there exists a $2$-cochain $m$ which assigns $1$ to every elementary domain that does not contain any $X$ or $O$ markings, such that the sign assignment $s'_2$ obtained by modifying $s_1$ by the $2$-cochain $m$ satisfies 
$h_{s_2}=h_{s'_2}$ and $v_{s_2}=v_{s'_2}$. Therefore, by Theorem \ref{thm:uniqueness}, $s'_2$ is gauge equivalent to $s_2$.
\end{proof}

\begin{thm}
The map $f$ from the set of orientation systems to the set of sign assignments induces a well-defined bijection $\wt{f}$ from the set of weak equivalence classes of orientation systems to the set of weak equivalence classes of sign assignments.
\end{thm}

\begin{proof}
Recall that two orientation systems $\mf{o}_1$ and $\mf{o}_2$ are weakly equivalent if and only if, for a fixed generator $x\in\mc{G}_{\mc{H}}$, $\mf{o}_1$ agrees with $\mf{o}_2$ on all the domains in $\mc{D}(x,x)$ that correspond to the empty periodic domains of $\mc{P}^0_{\mc{H}}$. Therefore, we need to find a basis for the empty periodic domains.

For each $i\in\{1,\ldots,l\}$, let $P_i=\sum_{j\in V(L_i)}V_j-\sum_{j\in H(L_i)}H_j$. These $l$ empty periodic domains  generate $\mc{P}^0_{\mc{H}}$, and $\sum_i P_i=0$ is the only relation among these domains. Therefore, the domains $P_1,\ldots,P_{l-1}$ freely generate $\mc{P}^0_{\mc{H}}$.

If $D\in\mc{D}(x,x)$ is a domain which corresponds to a vertical annulus $V_i$, then we know from Paragraph \ref{para:special} that $\mf{o}_1$ agrees with $\mf{o}_2$ on $D$ if and only if $v_{f(\mf{o}_1)}(i) = v_{f(\mf{o}_2)}(i)$. A similar statement holds for the horizontal annuli. A repeated application of the same principle shows that if $D\in\mc{D}(x,x)$ corresponds to the empty periodic domain $P_i$, then $\mf{o}_1$ agrees with $\mf{o}_2$ on $D$ if and only if $r_{f(\mf{o}_1)}(L_i)=r_{f(\mf{o}_1)}(L_i)$. Therefore, the orientation systems $\mf{o}_1$ and $\mf{o}_2$ are weakly equivalent if and only if the sign assignments $f(\mf{o}_1)$ and $f(\mf{o}_2)$ are weakly equivalent. This shows that the map in question is well-defined and injective. As both sets have $2^{l-1}$ elements, it is a bijection.
\end{proof}

A consequence of the theorems in this section is the following.

\begin{thm}
  There is a bijection $\wt{f}$ between the weak equivalence classes   of orientation systems and the weak equivalence classes of sign   assignments, such that for each of the $2^{l-1}$ weak equivalence   classes of orientation systems $\mf{o}$, the homology of the grid   chain complex, evaluated with the sign assignment $f(\mf{o})$, is   isomorphic as an absolutely $(l+1)$-graded group to   $\whhfl(L,\Z,\mf{o})\otimes_i(\otimes^{m_i-1}Q_i)$.
\end{thm}

\begin{figure}
\psfrag{x}{$X$}
\psfrag{o}{$O$}
\psfrag{xo}{$XO$}
\begin{center}
\includegraphics[width=0.8\textwidth]{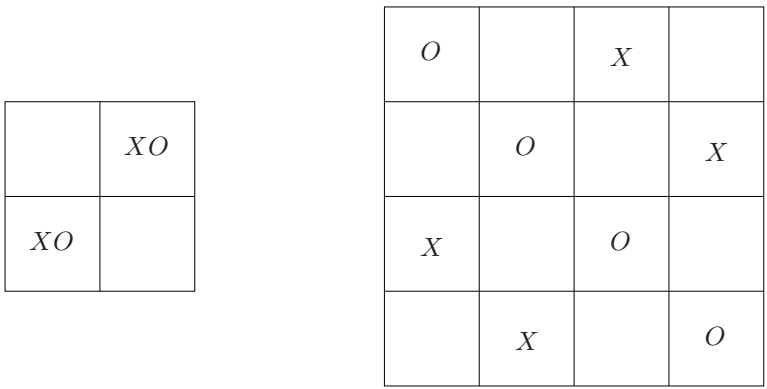}
\end{center}
\caption{Grid diagrams for the two-component unlink and the Hopf link.}\label{fig:grid}
\end{figure}

Let us conclude with a couple of examples. The first grid diagram in Figure \ref{fig:grid} represents the two-component unlink. There are exactly two generators and exactly two rectangles connecting the two generators. One weak equivalence class assigns the same sign to both the rectangles while the other weak equivalence class assigns opposite signs. Therefore, for one weak equivalence class of orientation systems, the homology is $\Z/2\Z$, while for the other other weak equivalence class of orientation systems, the homology is $\Z\oplus\Z$.

The second grid diagram in Figure \ref{fig:grid} represents the Hopf link. There are twenty-four generators and sixteen rectangles. It can be checked by direct computation that the homology is independent of the sign assignment. Therefore, the link Floer homology of the Hopf  link is the same for both the weak equivalence classes of orientation systems.

\bibliographystyle{amsplain}

\bibliography{signs}

\end{document}